\documentclass{notices}

\usepackage[T1]{fontenc}
\usepackage{amsmath}
\usepackage{amscd}
\usepackage{graphics}
\usepackage{latexsym}
\usepackage[all]{xy}
\usepackage{epigraph}
\usepackage{comment}
\usepackage{graphicx}
\usepackage{amsmath,amssymb,xy,array}
\usepackage[T1]{fontenc}
\usepackage{amsfonts}
\usepackage{amsmath}
\usepackage{amssymb}
\usepackage{amsthm}
\usepackage{graphicx}
\usepackage{amscd}
\usepackage{latexsym}
\usepackage{enumitem}
\usepackage{xcolor}
\usepackage{amsfonts}
\usepackage{xypic}
\usepackage[utf8]{inputenc}
\usepackage{hyperref}
\usepackage{mathrsfs}
\usepackage[toc,page]{appendix}
\usepackage{fancyhdr}
\usepackage{tikz}
\usepackage{tikz-cd}
\usepackage{longtable}
\usepackage{geometry}
\usepackage{textcomp}
\usetikzlibrary{trees}
\usetikzlibrary{matrix, arrows}
\usepackage{multirow}
\usepackage{youngtab}
\usepackage{mathdots}
\usepackage{tcolorbox}

\definecolor{yqyqyq}{rgb}{0.5019607843137255,0.5019607843137255,0.5019607843137255}\definecolor{uuuuuu}{rgb}{0.26666666666666666,0.26666666666666666,0.26666666666666666}
\definecolor{uququq}{rgb}{0.25098039215686274,0.25098039215686274,0.25098039215686274}
\definecolor{wwwwww}{rgb}{0.4,0.4,0.4}
\setlist[itemize]{leftmargin=6mm}

\renewcommand{\P}{\mathbb P}

\DeclareMathOperator{\NE}{NE}

\DeclareMathOperator{\Exc}{Exc}

\DeclareMathOperator{\MCD}{MCD}

\DeclareMathOperator{\Eff}{Eff}
\DeclareMathOperator{\Nef}{Nef}
\DeclareMathOperator{\Mov}{Mov}
\DeclareMathOperator{\Pic}{Pic}

\renewcommand{\sec}{\mathbb{S}ec}

\DeclareMathOperator{\Amp}{Amp}

\renewcommand{\P}{\mathbb{P}}

\newtheorem{thm}{Theorem}[section]

\theoremstyle{definition}

\newtheorem{Definition}[thm]{Definition}

\newtheorem{Example}[thm]{Example}

\hypersetup{pdfpagemode=UseNone}
\hypersetup{pdfstartview=FitH}

\title{On the weak Lefschetz principle\\ in birational geometry}

\date{}
\author{
  C\'esar Lozano Huerta
  \affil{
    CONACYT Research Fellow in mathematics at the National University of Mexico.
Email: lozano@im.unam.mx
    }
  \and
  Alex Massarenti
  \affil{
    Research Fellow in mathematics at the Department of Mathematics and Computer Science of the University of Ferrara (Italy).
His email address is alex.massarenti@unife.it
   }
}

\begin{document}

\maketitle

\section*{The lay of the land}

\noindent
Algebraic geometers study the solution sets to finite collections of polynomial equations in $\mathbb{C}^n$ or in the complex projective space $\mathbb{P}^n$. Often, the goal is to understand the basic geometric properties of these sets such as connectedness, smoothness, dimension, etc. The most important and familiar examples are solution sets to finite collections of linear equations and in this case, we notice that their geometry exhibits many similarities to that of their ambient spaces $\mathbb{C}^n$ or $\mathbb{P}^n$. 

\noindent
When one analyzes solution sets of higher degree equations, the situation changes radically. These sets no longer have an easily recognizable structure, such as that of a vector space or a group. Unlike the linear case, they may have many connected components and its dimension may be challenging to compute. Despite all these difficulties, this article will discuss scenarios in which the geometry of the solution set is similar to that of the space in which it sits. Such situations are somewhat unexpected and have had a profound impact in mathematics.

\noindent
Our departing point is the influential work of Solomon Lefschetz started in 1924 \cite{AS}. Consider a collection of finitely many homogeneous polynomials whose solution set in $\mathbb{P}^n$ has dimension bigger than one. Lefschetz analyzed the following: if one considers an additional generic linear equation and looks at the resulting solution set, the so-called hyperplane section, often one finds that many geometric properties of it are determined by the original solution set.

\noindent
In what follows, we will look at the hyperplane theorem of Lefschetz in different contexts. We will revisit its original formulation in algebraic topology, and build up to recent developments of it in birational geometry. In doing so, we will emphasize the main theme of this article: there are many contexts in geometry in which analogous versions of the Lefschetz hyperplane theorem hold. Even though it was originally stated in the context of algebraic topology, there are versions of it in algebraic geometry, homotopy theory, complex geometry and birational geometry. Intuitively speaking, all these results focus on a solution set (or a submanifold, depending on the context) whose geometric properties are determined by the space in which it sits.

\noindent
We have tried to present the development of the ideas somehow chronologically. We apologize for the omissions we unintentionally made due to our lack of historical understanding.

\section*{Notions from algebraic topology}
\noindent
Started in the early 20th century and known as combinatorial topology until the early 1940s, the field of algebraic topology constructs algebraic invariants in order to distinguish two topological spaces. Homology and cohomology groups are two such invariants and central notions of this article. These groups, whose definitions are subtle, turn out to be accessible to computations and are for this reason commonly used.
Let us recall some of their properties. 

\noindent
The $k$th homology group of a topological space $X$, with coefficients in a field $\mathbb{F}$, can be defined in terms of cycles and boundaries when a triangulation of $X$ is available. Indeed, the $k$th homology group of $X$ is defined as the quotient
$$H_k(X,\mathbb{F})=\{k\mbox{-cycles of }X\}/\{k\mbox{-boundaries of }X\}$$
where the coefficients in question are in the field $\mathbb{F}$. In this case, the group $H_k(X,\mathbb{F})$ has the structure of a vector space over $\mathbb{F}$ and in the remainder of this article our homology groups will have $\mathbb{Q}$-coefficients.

\noindent
A key feature of homology is the following: any continuous morphism between topological spaces $f:X\rightarrow Y$ induces a corresponding morphism 
$$f_*:H_k(X,\mathbb{Q})\rightarrow H_k(Y, \mathbb{Q})$$ 
which is a linear transformation between vector spaces over $\mathbb{Q}$.  

\noindent 
Alongside the $k$th homology group there is the $k$th cohomology group denoted by $H^k(X,\mathbb{Q})$. If $X$ is an oriented manifold, the former is the dual vector space to the latter, but the following property of cohomology is different: any continuous morphism between topological spaces $f:X\rightarrow Y$ induces a corresponding morphism
$$f^*:H^k(Y,\mathbb{Q})\rightarrow H^k(X, \mathbb{Q})$$
which is a linear transformation between vector spaces over $\mathbb{Q}$.

\noindent 
Cohomology has a distinctive advantage over homology: not only can we add cohomology classes, we can multiply them as well. When $X$ is an oriented manifold, this product has a rich geometric meaning which is often given by the intersection of two submanifolds. (We could do this in homology, but the labels in cohomology work better. For one thing, they yield a graded ring: if $\alpha \in H^k(X,\mathbb{Q}) $ and $\beta\in H^l(X,\mathbb{Q})$, then $\alpha.\beta\in H^{k+l}(X,\mathbb{Q}$)). This is the so-called cup product and will appear a bit later.

\section*{Tracing back the Lefschetz hyperplane theorem}
\noindent
In this section we formulate the theorem that sets the theme of this article: the Lefschetz hyperplane theorem. We start by describing its context.

\noindent
Let us consider a smooth complex algebraic surface $S\subset \mathbb{P}^n$; such an $S$ is a manifold of dimension $4$ over $\mathbb{R}$. We are interested in describing the homology group $H_1(S,\mathbb{Q})$ of the surface $S$ with coefficients in the rationals.

\noindent
Since the surface $S$ lies inside the projective space $\mathbb{P}^n$, one may take a smooth \textit{hyperplane section} of it by throwing in an additional linear equation to the defining equations of $S$ such that the result is a smooth algebraic curve $C\subset S$ of genus $g$. One can also consider a generic family of hyperplane sections of $S$ containing $C$ as follows. Linear polynomials in $\mathbb{P}^n$ form a family isomorphic to $\mathbb{P}^n$ for their $n+1$ coefficients can be thought of as coordinates. One point $p$ in this family gives rise to $C$. Pick another general point $q$ of this family and consider the unique line $\overline{pq}=L$ that contains them.  We now may think of the hyperplane sections of $S$ parametrized by $t\in L\cong \mathbb{P}^1$. That is, for each value $t\in L\cong\mathbb{P}^1$ there is a curve $C_t\subset S$.\footnote{The family $\{C_t\}$ is called a Lefschetz pencil and its key property is that its generic member is smooth and the singular members have all mild singularities.}

\noindent
The advantage of considering $C$ within this family is that one may try to describe the homology of the surface $H_1(S,\mathbb{Q})$ as a linear combination of the $1$-cycles of $C_t$. That is, we may consider a $1$-cycle of $C_t$ and then think of it as a $1$-cycle of $S$. A result by Lefschetz says this is possible: \textit{any $1$-cycle of $S$ is homologous on $S$ to a $1$-cycle lying on a generic $C_t$}. Since the homology of the generic $C_t$ is $H_1(C_t,\mathbb{Q})\cong \mathbb{Q}^{2g}$, this yields $H_1(S,\mathbb{Q})$ as a quotient of a well-known space. 
The following figure 
$$
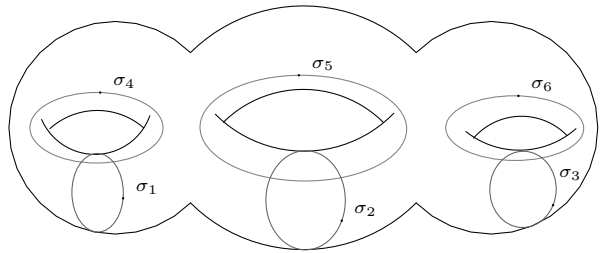
\definecolor{wqwqwq}{rgb}{0.3764705882352941,0.3764705882352941,0.3764705882352941}
\definecolor{yqyqyq}{rgb}{0.5019607843137255,0.5019607843137255,0.5019607843137255}
\begin{tikzpicture}[line cap=round,line join=round,>=triangle 45,x=1.0cm,y=1.0cm,scale=0.5]
\clip(1.,-0.4) rectangle (17.,6.4);
\draw [shift={(4.,3.)},line width=0.4pt]  plot[domain=0.7853981633974483:5.497787143782138,variable=\t]({1.*2.8284271247461903*cos(\t r)+0.*2.8284271247461903*sin(\t r)},{0.*2.8284271247461903*cos(\t r)+1.*2.8284271247461903*sin(\t r)});
\draw [shift={(9.,2.)},line width=0.4pt]  plot[domain=0.7853981633974483:2.356194490192345,variable=\t]({1.*4.242640687119286*cos(\t r)+0.*4.242640687119286*sin(\t r)},{0.*4.242640687119286*cos(\t r)+1.*4.242640687119286*sin(\t r)});
\draw [shift={(9.,4.)},line width=0.4pt]  plot[domain=3.9269908169872414:5.497787143782138,variable=\t]({1.*4.242640687119286*cos(\t r)+0.*4.242640687119286*sin(\t r)},{0.*4.242640687119286*cos(\t r)+1.*4.242640687119286*sin(\t r)});
\draw [shift={(14.,3.)},line width=0.4pt]  plot[domain=-2.356194490192345:2.356194490192345,variable=\t]({1.*2.8284271247461903*cos(\t r)+0.*2.8284271247461903*sin(\t r)},{0.*2.8284271247461903*cos(\t r)+1.*2.8284271247461903*sin(\t r)});
\draw [rotate around={0.:(3.5,3.)},line width=0.4pt,color=yqyqyq] (3.5,3.) ellipse (1.7695514389270113cm and 0.9387823469839274cm);
\draw [rotate around={-0.7282969530764222:(9.,2.99)},line width=0.4pt,color=yqyqyq] (9.,2.99) ellipse (2.747412709380753cm and 1.406334453701112cm);
\draw [rotate around={-0.3536731592464616:(14.62,2.99)},line width=0.4pt,color=yqyqyq] (14.62,2.99) ellipse (1.834112577519103cm and 0.8599238030278251cm);
\draw [shift={(3.46,3.84)},line width=0.4pt]  plot[domain=3.5532639970216877:5.940547237617583,variable=\t]({1.*1.549451515859725*cos(\t r)+0.*1.549451515859725*sin(\t r)},{0.*1.549451515859725*cos(\t r)+1.*1.549451515859725*sin(\t r)});
\draw [shift={(3.52,1.58)},line width=0.4pt]  plot[domain=0.8516984491879622:2.304957322788926,variable=\t]({1.*1.8804386558311053*cos(\t r)+0.*1.8804386558311053*sin(\t r)},{0.*1.8804386558311053*cos(\t r)+1.*1.8804386558311053*sin(\t r)});
\draw [shift={(9.02,5.7)},line width=0.4pt]  plot[domain=3.9227355235227654:5.513785778658482,variable=\t]({1.*3.3234319610908236*cos(\t r)+0.*3.3234319610908236*sin(\t r)},{0.*3.3234319610908236*cos(\t r)+1.*3.3234319610908236*sin(\t r)});
\draw [shift={(9.,1.)},line width=0.4pt]  plot[domain=0.7857052833686552:2.3519391967278684,variable=\t]({1.*3.02109501167155*cos(\t r)+0.*3.02109501167155*sin(\t r)},{0.*3.02109501167155*cos(\t r)+1.*3.02109501167155*sin(\t r)});
\draw [shift={(14.74,4.68)},line width=0.4pt]  plot[domain=4.039019779413229:5.43688693550415,variable=\t]({1.*2.277015590636129*cos(\t r)+0.*2.277015590636129*sin(\t r)},{0.*2.277015590636129*cos(\t r)+1.*2.277015590636129*sin(\t r)});
\draw [shift={(14.8,1.62)},line width=0.4pt]  plot[domain=0.7682999593860934:2.42853182073094,variable=\t]({1.*1.6861767457329406*cos(\t r)+0.*1.6861767457329406*sin(\t r)},{0.*1.6861767457329406*cos(\t r)+1.*1.6861767457329406*sin(\t r)});
\draw [rotate around={-89.26228927654937:(3.529819596244108,1.2706377359139842)},line width=0.4pt,color=wqwqwq] (3.529819596244108,1.2706377359139842) ellipse (1.0450082197728656cm and 0.6832499610838083cm);
\draw [rotate around={90.:(9.06,1.07)},line width=0.4pt,color=wqwqwq] (9.06,1.07) ellipse (1.31836477364618cm and 1.0554551986660239cm);
\draw [rotate around={87.70938995736422:(14.84,1.36)},line width=0.4pt,color=wqwqwq] (14.84,1.36) ellipse (1.0141884227859228cm and 0.8821440681164167cm);
\begin{scriptsize}
\draw [color=black] (4.207320059432352,1.1228127668908174) circle (0.5pt);
\draw[color=black] (4.8503246716281065,1.3873927893843265) node {$\sigma_1$};
\draw [color=black] (10.021646014634038,0.5266466970209462) circle (0.5pt);
\draw[color=black] (10.655280588890975,0.8027210206784856) node {$\sigma_2$};
\draw [color=black] (15.639127791391797,0.9403535114528692) circle (0.5pt);
\draw[color=black] (16.1,1.75) node {$\sigma_3$};
\draw [color=black] (3.599672650523572,3.9372919381036433) circle (0.5pt);
\draw[color=black] (4.244771716230181,4.206345959930345) node {$\sigma_4$};
\draw [color=black] (8.880564751751765,4.396208977009445) circle (0.5pt);
\draw[color=black] (9.527699223667252,4.665730921056363) node {$\sigma_5$};
\draw [color=black] (14.719040692669223,3.8482049112336933) circle (0.5pt);
\draw[color=black] (15.353536277323153,4.1228214215437955) node {$\sigma_6$};
\end{scriptsize}
\end{tikzpicture}
$$
depicts the generating $1$-cycles of $H_1(C,\mathbb{Z})$, where $C$ is a smooth curve of genus $3$.

\noindent
Let us state again the previous result. Let $i:C\hookrightarrow S$ be the inclusion of a smooth hyperplane section of $S$. Then the result above asserts that the induced map $$i_*:H_1(C,\mathbb{Q})\rightarrow H_1(S,\mathbb{Q})$$ is surjective. In this setting, this is the Lefschetz hyperplane theorem. In order to fully describe the first homology group of $S$ we need to further investigate the kernel of $i_*$; which consists of the so-called \textit{vanishing cycles}. The necessary information about such cycles can now be extracted from thinking of $C$ in the family $\{C_t\}$. From this analysis, one gets an important theorem: the dimension of $H_1(S,\mathbb{Q})$, called the first Betti number of $S$ and denoted by $b_1(S)$, is even.

\noindent
In general it is quite challenging to deduce properties of $1$-cycles on $S$ from their behavior as $1$-cycles on the curves $C_t$.
In fact, Zariski points out that such arguments are the most difficult part of all the topological arguments by Lefschetz \cite[p. 136]{ZARISKI}. For example, Zariski comments that Lefschetz's topological proof of the fact that the Betti number $b_1(S)$ is even seems to be incomplete.\footnote{Zariski comments that Lefschetz' proof is based on an unproven statement of Picard.}

\noindent
This difficulty, however, can be overcome via Hodge theory in which case the key is the following claim: intersection with the curve $C_t$ defines an isomorphism
\stepcounter{thm}
\begin{equation}\label{iso1}
\begin{aligned}
\tilde{\phi} &: H_3(S,\mathbb{Q})\rightarrow H_1(S,\mathbb{Q}).\\
\end{aligned}
\end{equation}

\noindent
The proof that $\tilde{\phi}$ is an isomorphism was first completed in \cite{HODGE}. We will come back to this map in the next section.

\noindent
Notice that $\tilde{\phi}$ factors through $H_1(C_t,\mathbb{Q})$. 
In fact, one can deduce that $\tilde{\phi}$ is an isomorphism from the surjectivity of $i_*$
\stepcounter{thm}
\begin{equation}\label{iso2}
\begin{aligned}
i_* &: H_1(C_t,\mathbb{Q})\rightarrow H_1(S,\mathbb{Q})\\
\end{aligned}
\end{equation}
by verifying that its kernel consists of vanishing cycles and is orthogonal to the cycles that generate $H_3(S,\mathbb{Q})$.

\noindent
The previous discussion is part of Lefschetz's original approach to the study of the topology of an algebraic surface $S$.
In what follows we will examine a generalization of this discussion which goes under the name of  the weak Lefschetz theorem. Before that, we decided to include a brief description of the strong Lefschetz theorem so the reader may contrast. Strong Lefschetz will concern the map $\tilde{\phi}$ while weak Lefschetz will concern the map $i_*$ above. 

\section*{Strong Lefschetz theorem}
\noindent 
The previous isomorphism (\ref{iso1}) admits a dual formulation involving cohomology which builds up to the Strong Lefschetz Theorem. Indeed, the cup product with the cohomology class $\omega$ of the hyperplane section $C_t$ yields the map $$\phi:H^1(S,\mathbb{Q})\rightarrow H^3(S,\mathbb{Q})$$
and $\phi$ is an isomorphism.  A generalization now may be formulated as follows. Let $X$ be a smooth complex algebraic variety of  dimension $n$. Then the map defined by the cup product with the powers $\omega^{n-k}$ 
$$\phi^{n-k}: H^k(X,\mathbb{Q})\rightarrow H^{2n-k}(X,\mathbb{Q})$$ 
is an isomorphism \cite{HODGE}. This result by Hodge is carried out over the field of complex numbers. In characteristic $p>0$ 
it was studied by Grothendieck, Artin and Verdier in \cite{SGA4} and they called it the strong\footnote{This is also known as the `hard Lefschetz theorem'.} Lefschetz theorem.

\section*{Weak Lefschetz theorem}\label{weakLefschetz}

\noindent
Let us keep using cohomology notation in this section so the contrast between the weak and strong versions of the Lefschetz theorem becomes clear. The weak Lefschetz theorem, in our context, claims that if $i : H \hookrightarrow X$ is the inclusion morphism of a smooth hyperplane section, then the induced map on the $k$th cohomology groups
\stepcounter{thm}
\begin{equation}\label{iso3}
\begin{aligned}
i^* &: H^k(X, \mathbb{Q}) \rightarrow H^k(H,\mathbb{Q})\\
\end{aligned}
\end{equation}
is an isomorphism for $k \le n- 2$ and an injection for $k = n - 1$. For complex varieties, the theorem follows from the exact sequence 
$$H^k_c(U, \mathbb{Q})\rightarrow H^k(X, \mathbb{Q})\rightarrow H^k(H, \mathbb{Q})\rightarrow H^{k+1}_c(U, \mathbb{Q})$$ 
where $U$ denotes the open set $U=X\setminus H$ and $H^*_c(U, \mathbb{Q})$ stands for cohomology with compact support. Since $H_c^k(U,\mathbb{Q}) \cong H_{2n-k}(U,\mathbb{Q})$, the weak Lefschetz theorem is therefore equivalent to the vanishing of the homology groups: $H_j(U, \mathbb{Q}) =0$ for $j\ge n + 1$.

\noindent
The analogous result in characteristic $p>0$ was proved by Grothendieck, Artin and Verdier \cite{SGA4} in the context of \'etale cohomology with $\mathbb{Q}_l$-adic coefficients using exact sequences and cohomology vanishing.
They called it the \textit{weak Lefschetz theorem}.

\noindent
It is important to mention that the Lefschetz theorems, as formulated above, bear a close relation to the Weil conjectures. We refer the reader to \cite{MR241437}, in particular Kleiman's \textit{expos\'e X},  for an account about this fascinating part of algebraic geometry.

\noindent
The isomorphism (\ref{iso3}) is a crucial ingredient in what follows. In fact, we will exhibit examples where this isomorphism also preserves objects contained in the $2$nd cohomology group which encode birational information about the varieties $X$ and $Y$. This is the `weak Lefschetz' of our title. In the following section we will motivate the word `principle'.

\section*{Building up to a principle}\label{principle}
\noindent
A different approach to study the topology of a smooth projective variety $X$ was developed by Bott \cite{BOTT} and by Andreotti-Frankel \cite{ANDRE} both in 1959. These two papers used Morse theory applied to a non-degenerate real $C^{\infty}$-function $f$ which satisfies $f(x)=0$ if $x \in H$ and $f(x)>0$ if $x\notin H$, where $H\subset X$ is a smooth hyperplane section. Bott even describes a cellular structure on $X$ by attaching cells to $H$. This decomposition allows him to analyze the topology of $X$ in detail. In particular, Bott proves the statement dual to the weak Lefschetz theorem above: the induced morphism by the inclusion $$i_*:H_k(H,\mathbb{Q})\rightarrow H_k(X,\mathbb{Q})$$ is bijective if $0\le k \le n-2$, and $$i_*: H_{n-1}(H,\mathbb{Q})\rightarrow H_{n-1}(X,\mathbb{Q})$$ is surjective. 

\noindent
Bott actually proves more. He argues that one may focus on other topological invariants instead of homology or cohomology. 
For example, Bott applied this theory to homotopy groups and showed that the inclusion $i:H\hookrightarrow X$ yields isomorphisms for the homotopy groups 
$$i_*:\pi_k(H)\rightarrow \pi_k(X).$$ 
This result, along with those for homology, starts looking more like a \textit{principle}. Indeed, the relationship between the geometry of a variety $X$ and that of a hyperplane section $H$ seems to be robust as we are planning to exhibit in the rest of this article.

\section*{The weak Lefschetz principle}\label{WLP}
\noindent
After the previous sections we have the following situation. Let $X, Y$ be smooth complex projective varieties equipped with embedding $i: Y\hookrightarrow X$. This inclusion $i$ may induce an isomorphism on the cohomology, homology or homotopy groups. It turns out that it may also induce an isomorphism on objects in some other contexts, for example the Picard groups, \'etale cohomology or Hodge structures just to name a few. One thus may be led to think that the weak Lefschetz theorem seems to be robust, to the extent of considering it part of a principle. This principle says that often the geometry of the ambient variety $X$ determines that of the subvariety $Y$. Finding the precise conditions for this principle to hold is subject of current research. 

\noindent
In the rest of the article, we aim to exhibit examples in the context of birational geometry for which the weak Lefschetz principle holds. The next section contains the precise definitions to this end and is the technical part of the article. It may be skipped on a first reading.

\section*{Notions from birational geometry}\label{bir_geo}
Now we change the gears from topology to birational geometry; that is, we return to studying zeros of polynomials. The set of all algebraic varieties is vast, so one tries to organize it into equivalence classes. Birational geometry studies algebraic varieties up to birational equivalence and often seeks to find the ``simpliest'' representative in each equivalence class.

\noindent
From now on, we assume some familiarity with elementary concepts from algebraic geometry such as Cartier divisor, linear equivalence, blowup and strict transform. Since the definitions introduced in this section might be hard to grasp, we have included an explicit example that can be read independently in Boxes $1$, $2$, and $3$.

\noindent
Let us consider $X$ a smooth projective variety over $\mathbb{C}$. If $D$ is a subvariety of codimension $1$ in $X$, then one can naturally think of the number $D \cdot C$, the result of intersecting $D$ with a curve $C$. Loosely speaking, should count the number of points in the intersection between $D$ and $C$. However, if $C$ is a subset of $D$ we need a different description. Let us recall its formal definition so we are on firm ground. Given a Cartier divisor $D$ on a normal variety $X$ and a curve $C\subset X$ their intersection number is defined as
\stepcounter{thm}
\begin{equation}\label{int}
D\cdot C = \deg\mathcal{O}_C(D),
\end{equation}
where $\mathcal{O}_C(D)$ is the restriction to $C$ of the invertible sheaf $\mathcal{O}_X(D)$ associated to $D$. When $D$ does not contain $C$ the number $D\cdot C$ is the number of points of intersection between $D$ and $C$ counted with multiplicity. Note that this definition can be extended by linearity to $\mathbb{Q}$-Cartier divisors, that is Weil divisors having a multiple that is Cartier, and arbitrary linear combinations of irreducible curves; which are called $1$-cycles of $X$.

\noindent
Now the correct notion of intersection reveals structure in the set of divisors. Indeed, we say that two Cartier divisors $D,D'\subset X$ are numerically equivalent if $D\cdot C = D'\cdot C$, for any irreducible curve $C\subset X$. In this case, we write $D\equiv D'$. The group of all Cartier divisors on $X$ modulo numerical equivalence is a finitely generated free abelian group, called the Neron-Severi group and it is denoted by $N^1(X)$. Hence, the vector space $$N^1(X)_{\mathbb{R}} = N^1(X)\otimes\mathbb{R}$$ is finite dimensional \cite[page 8]{De01} and its dimension, denoted by $\rho_X$, is called the \textit{Picard rank of $X$}.

\noindent
The study of the Picard rank has deep roots in algebraic geometry. When $X$ is a smooth surface over $\mathbb{C}$, the number $\rho_X$ was carefully analyzed by Picard, Severi and Lefschetz. In fact, Lefschetz observed that curves in $X $, up to numerical equivalence, can be thought of as elements in the cohomology group $H^2(X,\mathbb{Q})$. From this he was able to show that $\rho_X$  is less or equal than the dimension of $H^2(X,\mathbb{Q})$; hence finite. This enhances the analysis mentioned earlier in this article about the first Betti number $b_1(S)$ of an algebraic surface $S$.

\noindent
So far we have mentioned three equivalence relations that can be imposed on Cartier divisors: linearly equivalent, numerically equivalent, and cohomologous. These three relations may coincide, for instance if $X=\mathbb{P}^n$, but in general they are different; for example on a general quartic surface in $\mathbb{P}^3$. The first two relations give rise to the Picard group $\Pic(X)$ and $N^1(X)$, respectively. Later in this section, we will restrict our study to spaces where linear and numerical equivalence coincide.

\noindent
Let us now recall the definition of some cones contained in the space $N^1(X)_{\mathbb{R}}$ which are important in birational geometry and central for the rest of the article. See the Boxes $1$, $2$, and $3$ for an example.

\begin{Definition}
Let $D\subset X$ be a $\mathbb{Q}$-Cartier divisor, that is a divisor such that $mD$ is Cartier for some integer $m>0$.
\begin{itemize}
\item[-] The divisor class is said to be \textit{effective} if it represents an actual subvariety of codimension $1$. The \emph{effective cone} of $X$ is the convex cone $\Eff(X)\subset N^1(X)_{\mathbb{R}}$ generated by classes of effective divisors. 
\item[-] $D$ is \textit{very ample} if it induces an embedding, and $D$ is \textit{ample} if $mD$ is very ample for $m\gg 0$. The \emph{ample cone} of $X$ is the convex cone $\Amp(X)\subset N^1(X)_{\mathbb{R}}$ generated by classes of ample divisors. 
\item[-] $D$ is \textit{nef} if $D\cdot C\geq 0$ for any irreducible curve $C\subset X$. Note that if $D$ is ample then $D\cdot C > 0$ for any curve $C\subset X$. Hence nefness is a mild relaxation of ampleness. The \emph{nef cone} of $X$ is the closed convex cone $\Nef(X)\subset N^1(X)_{\mathbb{R}}$ generated by classes of nef divisors. 
\item[-] The stable base locus of $D$ is the set of points $p\in X$ such that for all $m > 0$, if $mD$ is integral, all the divisors in the linear system of $mD$ pass through $p$. The \emph{movable cone} of $X$ is the convex cone $\Mov(X)$ generated by classes of \emph{movable divisors}; these are divisors whose stable base locus has codimension at least two in $X$.
\end{itemize}
\end{Definition}

\noindent
We have inclusions among the previous cones
$$\Amp(X)\subset\Nef(X)\subset \overline{\Mov(X)} \subset \overline{\Eff(X)}\subset N^1(X)_{\mathbb{R}}.$$

\medskip
\begin{tcolorbox}
\begin{center}\textbf{Box 1: the effective cone and the Mori cone of curves}\end{center}
Let us work out explicitly the cone of effective divisors and the Mori cone of curves of $X$, the blow-up of $\mathbb{P}^n$ at two points $p,q\in\mathbb{P}^n$, with $n>1$. Let $H,H_p,H_q,H_{p,q}$ be the strict transforms respectively of a hyperplane, a hyperplane passing through $p$, through $q$, and through both $p$ and $q$. Moreover, let $E_p, E_q$ be the exceptional divisors over $p$ and $q$ respectively. Note that $H_p = H-E_p, H_q = H-E_q$ and $H_{p,q} = H-E_p-E_q$. Then $N^1(X)\cong \mathbb{Z}[H,E_p,E_q]$. We will denote by $h$ the strict transform of a general line in $\mathbb{P}^n$, and by $e_p,e_q$ classes of lines in $E_p$ and $E_q$ respectively. The intersection pairing (\ref{intpai}) is given by $H\cdot h = 1, H\cdot e_p = H\cdot e_q = 0, E_p\cdot e_q = E_q\cdot e_p = 0, E_p\cdot e_p = E_q\cdot e_q = -1$. The last two intersections numbers might be not obvious from a geometrical point of view. To compute them one may reason as follows: the divisor $H-E_p$ represents the strict transform of a general hyperplane through $p$, and $h-e_p$ represents the strict transform of a general line through $p$. In the blow-up $X$ these strict transforms do not intersect anymore, so $0 = (H-E_p)\cdot (h-e_p) = H\cdot h-H\cdot e_p - E_p\cdot h + E_p\cdot e_p$ and hence $E_p\cdot e_p = - H\cdot h = -1$.

Now, let $C\subset X$ be an irreducible curve. Then either $C$ is contained in an exceptional divisor and then it is numerically equivalent to a positive multiple of $e_p$ or $e_q$, or it is mapped by the blow-down map to an irreducible curve $\Gamma\subset\mathbb{P}^n$. Let $d,m_p,m_q$ be respectively the degree and the multiplicities of $\Gamma$ at $p$ and $q$. Then $C\equiv dh-m_pe_p-m_qe_q$. We may write $C\equiv d(h-e_p-e_q)+(d-m_p)e_p+(d-m_q)e_q$. Furthermore, $d-m_p >0$ otherwise by B\'ezout's theorem $\Gamma$ would contain a line through $p$ as a component, and similarly $d-m_q > 0$. Hence $\NE(X)$ is closed and generated by the classes $e_p,e_q$ and $h-e_p-e_q$. Note that the latter is the strict transform of the line in $\mathbb{P}^n$ through $p,q$. Similarly, it can be shown that $\Eff(X)$ is closed and generated by the classes of $E_p, E_q$ and $H_{p,q}$.
\end{tcolorbox}

\newpage

\noindent
Similar to the case of divisors, two $1$-cycles $C,C'\subset X$ are numerically equivalent if $D\cdot C = D\cdot C'$ for any irreducible Cartier divisor $D\subset X$. We will denote by $N_1(X)$ the quotient of the group of $1$-cycles by this equivalent relation, and consider the vector space $N_1(X)_{\mathbb{R}} = N_1(X)\otimes\mathbb{R}$.

\noindent
Intersecting curves and divisors according to (\ref{int}) induces a non-degenerate pairing
\stepcounter{thm}
\begin{equation}\label{intpai}
N^1(X)\times N_1(X)\rightarrow\mathbb{Z}
\end{equation}
which implies that $N_1(X)_{\mathbb{R}}$ is finite dimensional.

\noindent
The following cone is called the Mori cone of curves and was introduced by S. Mori. Let $\mathrm{NE}(X)$ denote the closure of the convex cone in $N_1(X)_{\mathbb{R}}$ of classes of effective $1$-cycles, that is classes in $N^1(X)_{\mathbb{R}}$ which represent actual curves in $X$.

\bigskip

\begin{tcolorbox}
\begin{center}\textbf{Box 2: the Nef cone and movable cone (Box 1 continued)}\end{center}
Let us work out the nef cone of $X$ when $n>2$. This is the cone of divisors intersecting non negatively all the irreducible curves in $X$. 
Since any curve in $X$ can be written as a linear combination with non-negative coefficients of the generators of $\NE(X)$, it is enough to check when a divisor intersects non-negatively these generators. Let us write $D\equiv aH+bE_p+cE_q$. Then $D\cdot (h-e_p-e_q) = a+b+c$, $D\cdot e_p = -b$ and $D\cdot e_p = -c$, and $\Nef(X)$ is defined in $N^1(X)_{\mathbb{R}}\cong\mathbb{R}^3$ by the inequalities $a+b+c\geq 0, b\leq 0, c\leq 0$. Hence $\Nef(X)$ is generated by $\langle H, H_p, H_q\rangle$.

Finally, we determine the movable cone of $X$. The divisor $H_{p,q}$ represents the hyperplanes of $\mathbb{P}^n$ passing through $p,q$. Hence the stable base locus of $H_{p,q}$ consists of the strict transform of the line through $p,q$. The stable base locus of all divisors in the cone generated by $\langle H_p,H_q,H_{p,q}\rangle$ is contained in such a strict transform as $H_p,H_q$ have no base loci. Hence all the divisors in this cone are movable when $n>2$. On the other hand, all divisors in the interior of the cone $\langle H, H_p, E_q\rangle$ contain $E_q$, all divisors in the interior of the cone $\langle H, H_q, E_p\rangle$ contain $E_p$, and all divisors in the interior of the cone $\langle H, E_p, E_q\rangle$ contain $E_p\cup E_q$. Therefore, $\Mov(X)$ is the cone generated by $\langle H, H_p, H_q, H_{p,q}\rangle$. 
\end{tcolorbox}

\bigskip

\noindent
We now introduce the main notion of this section: Mori dream spaces in the context of the minimal model program. The goal of the minimal model program is to construct a birational model of any complex projective variety which is as simple as possible in a suitable sense. This subject has its origins in the classical birational geometry of surfaces studied by the Italian school. In $1988$ S. Mori extended the concept of minimal model to $3$-folds by allowing suitable singularities on them \cite{Mo88} and was awarded the Fields Medal for his contributions on this topic. In $2010$ there was another breakthrough in minimal model theory when C. Birkar, P. Cascini, C. Hacon and J. McKernan proved the existence of minimal models for a suitable class of higher dimensional varieties \cite{BCHM10}.
C. Birkar was awarded the Fields Medal in 2018 for his contributions to birational geometry.

\medskip\noindent
\textit{Mori dream spaces}, introduced by Y. Hu and S. Keel in $2000$ \cite{HK00}, form a class of algebraic varieties that behave very well from the point of view of the minimal model program. In order to recall their definition, let us first define $\mathbb{Q}$-factorial modifications.

\medskip\noindent
We say that a birational map  $f: X \dasharrow X'$ to a normal projective variety $X'$  is a \emph{birational contraction} if its inverse does not contract any divisor. 
We say that it is a \emph{small $\mathbb{Q}$-factorial modification} 
if $X'$ is $\mathbb{Q}$-factorial, that is any Weil divisor on $X'$ has a multiple which is Cartier, and $f$ is an isomorphism in codimension one. 
If  $f: X \dasharrow X'$ is a small $\mathbb{Q}$-factorial modification, then 
there is pull-back map $$f^*:N^1(X')\to N^1(X)$$ which sends $\Mov(X')$ and $\Eff(X')$
isomorphically onto $\Mov(X)$ and $\Eff(X)$, respectively. In particular, we have $f^*(\Nef(X'))\subset \overline{\Mov(X})$.

\noindent
The previous paragraph makes explicit the importance in birational geometry of the cones defined earlier in this section. The next definition will restrict our study to spaces whose birational geometry is encoded in finitely many of such cones.

\begin{Definition}\label{def:MDS}\cite[Definition 1.10]{HK00}
A normal projective $\mathbb{Q}$-factorial variety $X$ is called a \emph{Mori dream space}
if the following conditions hold:
\begin{enumerate}
\item[-] $\Pic(X)\otimes \mathbb{Q}=N^1(X)\otimes \mathbb{Q}$,
\item[-] $\Nef{(X)}$ is generated by the classes of finitely many semi-ample divisors, that is, divisors having a multiple with empty base locus,
\item[-] there is a finite collection of small $\mathbb{Q}$-factorial modifications
 $f_i: X \dasharrow X_i$ such that each $X_i$ satisfies the second condition above and $\Mov{(X)} \ = \ \bigcup_i \  f_i^*(\Nef{(X_i)})$.
\end{enumerate}
\end{Definition}

\noindent
The collection of all faces of all cones $f_i^*(\Nef{(X_i)})$ in Definition \ref{def:MDS} forms a wall-and-chamber decomposition of $\Mov(X)$.
If two maximal cones of this fan, say $f_i^*(\Nef{(X_i)})$ and $f_j^*(\Nef{(X_j)})$, meet along a facet,
then there exist a normal projective variety $Y$, small birational morphisms $h_i:X_i\rightarrow Y$ and $h_j:X_j\rightarrow Y$ of relative Picard rank one, and a small modification $\varphi:X_i\dasharrow X_j$ such that $h_j\circ\varphi = h_i$. The fan structure on $\Mov(X)$ can be extended to a fan supported on $\Eff(X)$ as follows. We refer the reader to \cite[Proposition 1.11]{HK00} and \cite[Section 2.2]{Ok16} for details and Box 3 for an example.

\begin{Definition}\label{MCD}
Let $X$ be a Mori dream space. A wall-and-chamber decomposition of the effective cone $\Eff(X)$, called the \emph{Mori chamber decomposition} and denoted $\MCD(X)$, is described as follows.
There are finitely many birational contractions from $X$ to Mori dream spaces, denoted by $g_i:X\dasharrow Y_i$.
The maximal cones $\mathcal{C}$ of the $\MCD(X)$ are of the form: $\mathcal{C}_i \ = g_i^*\big(\Nef(Y_i)\big) * \mathbb{R}_{\ge 0}\Exc(g_i)$, where $\mathrm{Exc}(g_i)$ is the exceptional locus of $g_i$. Here $A * B$ denotes the join of the cones $A$ and $B$. We call $\mathcal{C}_i$ a \emph{maximal chamber} of $\Eff(X)$. We thus have $\Eff(X)=\bigcup_j \mathcal{C}_j$.
\end{Definition}

\noindent
The importance of the previous definition is that it tells us precisely how the cones in $N^1(X)_{\mathbb{R}}$ defined earlier encode birational information of a Mori dream space $X$. For one thing, it tells us the behavior of such cones when $X$ undergoes a $\mathbb{Q}$-factorial modification.

\begin{tcolorbox}
\begin{center}\textbf{Box 3: the Mori chamber decomposition}\end{center}

We worked out explicitly the Mori chamber decomposition of the blow-up of $\mathbb{P}^n$ at two points $X = Bl_{p,q}\mathbb{P}^n$ in Boxes $1$ and $2$. The following picture is a two dimensional cross-section of $\Eff(X)$ displaying its Mori chamber decomposition:
$$
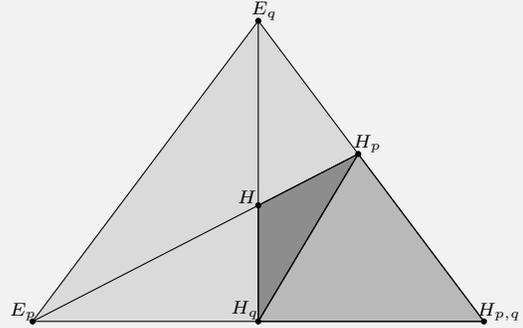
\begin{tikzpicture}[thick, scale=0.50][line cap=round,line join=round,>=triangle 45,x=1.0cm,y=1.0cm]

\fill[line width=0.4pt,fill=black,fill opacity=0.10000000149011612] (2.,2.) -- (14.,2.) -- (8.,10.) -- cycle;
\fill[line width=0.4pt,fill=black,fill opacity=0.3499999940395355] (8.,5.088435887573818) -- (10.657668310783288,6.456442252288949) -- (8.,2.) -- cycle;
\fill[line width=0.4pt,fill=black,fill opacity=0.15000000596046448] (8.,2.) -- (10.657668310783288,6.456442252288949) -- (14.,2.) -- cycle;
\draw [line width=0.4pt] (2.,2.)-- (14.,2.);
\draw [line width=0.4pt] (14.,2.)-- (8.,10.);
\draw [line width=0.4pt] (8.,10.)-- (2.,2.);
\draw [line width=0.4pt] (2.,2.)-- (10.657668310783288,6.45644225228895);
\draw [line width=0.4pt] (8.,10.)-- (8.,2.);
\draw [line width=0.4pt] (10.657668310783288,6.45644225228895)-- (8.,2.);
\draw [line width=0.4pt] (8.,5.088435887573818)-- (10.657668310783288,6.456442252288949);
\draw [line width=0.4pt] (10.657668310783288,6.456442252288949)-- (8.,2.);
\draw [line width=0.4pt] (8.,2.)-- (8.,5.088435887573818);
\draw [line width=0.4pt] (8.,2.)-- (10.657668310783288,6.456442252288949);
\draw [line width=0.4pt] (10.657668310783288,6.456442252288949)-- (14.,2.);
\draw [line width=0.4pt] (14.,2.)-- (8.,2.);
\begin{scriptsize}
\draw [fill=black] (2.,2.) circle (1.5pt);
\draw[color=black] (1.75,2.257731958762886) node {$E_{p}$};
\draw [fill=black] (14.,2.) circle (1.5pt);
\draw[color=black] (14.4,2.257731958762886) node {$H_{p,q}$};
\draw [fill=black] (8.,10.) circle (1.5pt);
\draw[color=black] (8.161258278145695,10.268041237113401) node {$E_{q}$};
\draw [fill=black] (8.,5.088435887573818) circle (1.5pt);
\draw[color=black] (7.7,5.3) node {$H$};
\draw [fill=black] (10.657668310783288,6.456442252288949) circle (1.5pt);
\draw[color=black] (10.9,6.72680412371134) node {$H_{p}$};
\draw [fill=black] (8.,2.) circle (1.5pt);
\draw[color=black] (7.65,2.3) node {$H_{q}$};
\end{scriptsize}
\end{tikzpicture}
$$
The divisors $H, H_p, H_q, H_{p,q}$ generate $\Mov(X)$, and $H, H_p, H_q$ generate $\Nef(X)$. The chamber delimited by $H,H_q,E_p$ corresponds to the contraction of $E_p$, similarly the chamber delimited by $H,H_p,E_q$ corresponds to the contraction of $E_q$, and chamber delimited by $H,E_p,E_q$ corresponds to the contraction of both $E_p$ and $E_q$.  

In the case $n\geq 3$, then $X$ admits only one small $\mathbb{Q}$-factorial modification $X'$ corresponding to the chamber delimited by $H_p,H_q,H_{p,q}$. In what follows (on this and the following page), we will investigate the geometry of $X'$. 

Consider a divisor lying on the wall delimited by $H_p,H_q$, for instance $D = H_p+H_q = 2H-E_p-E_q$ and let $L$ be the strict transform of the line through $p$ and $q$. Then $D\cdot L = 0$ and the linear system of quadrics in $\mathbb{P}^n$ through $p$ and $q$ induces a morphism $h_D:X\rightarrow Y$ contracting $L$ to a point. 

On the other hand, a divisor in the maximal chamber delimited by $H_p,H_q,H_{p,q}$ must be ample on $X'$. We can write such a divisor as $aH_p+bH_q+cH_{p,q}$ with $a,b,c > 0$ and observe that $(aH_p+bH_q+cH_{p,q})\cdot L = -c <0$. 

Note that the curve $L$ prevents divisors in the chamber $\left\langle H_p,H_q,H_{p,q}\right\rangle$ from being ample.

Let $g:W\rightarrow X$ be the blow-up of $X$ along $L$ with exceptional divisor $E_L\subset W$. Observe that $E_L$ is a $\mathbb{P}^{n-2}$-bundle over $L$. 
\end{tcolorbox}

\begin{tcolorbox}
There is a morphism $g':W\rightarrow X'$ contracting $E_L$, in the direction of $L$, onto a subvariety $Z\subset X'$ such that $Z\cong\mathbb{P}^{n-2}$. Consider the divisor $D' \equiv H_p+H_q+H_{p,q} \equiv 3H-2E_p-2E_q$. The linear system of $D'$ induces a rational map $\phi_{D'}:X\dasharrow X'$, and we have the following commutative diagram
 \[
  \begin{tikzpicture}[xscale=1.5,yscale=-1.2]
 \node (A1_2) at (1, -1) {$W$};
    \node (A0_0) at (0, 0) {$X$};
    \node (A0_2) at (2, 0) {$X'$};
    \node (A1_1) at (1, 1) {$Y$};
    \path (A0_0) edge [->,swap]node [auto] {$\scriptstyle{h_D}$} (A1_1);
    \path (A0_2) edge [->]node [auto] {$\scriptstyle{h}$} (A1_1);
    \path (A0_0) edge [->,dashed]node [auto] {$\scriptstyle{\varphi_{D'}}$}  (A0_2);
    \path (A1_2) edge [->,swap]node [auto] {$\scriptstyle{g}$}  (A0_0);
    \path (A1_2) edge [->]node [auto] {$\scriptstyle{g'}$}  (A0_2);
  \end{tikzpicture}
  \]    
where $h:X'\rightarrow Y$ is a small modification contracting $Z\subset X'$ to $h_D(L)$. The rational map $\phi_{D'}:X\dasharrow X'$ is an isomorphism between $X\setminus L$ and $X'\setminus Z$ and replaces $L$ with the variety $Z$ which is covered by curves having non-negative intersection with all divisors in the chamber $\left\langle H_p,H_q,H_{p,q}\right\rangle$. Concretely, in the case $n = 3$ for instance, we can fix homogeneous coordinates $[x:y:z:w]$ on $\mathbb{P}^3$, assume that $p = [1:0:0:0]$, $q = [0:0:0:1]$, and consider the rational maps
$$\alpha: \mathbb{P}^3\dasharrow \mathbb{P}^7$$
defined by $\alpha([x:y:z:w]) = [xy:xz:xw:y^2:yz:yw:z^2:zw]$, that is induced by the quadrics of $\mathbb{P}^{3}$ passing through $p$ and $q$, and 
$$\beta: \mathbb{P}^3\dasharrow \mathbb{P}^{11}$$ 
defined by $\beta([x:y:z:w]) = [xy^2:xz^2:xyz:xyw:xzw:y^3:y^2z:y^2w:yz^2:yzw:z^3,z^2w]$, that is induced by the cubics of $\mathbb{P}^{3}$ having at least double points at $p$ and $q$. Then $Y$ is the closure of the image of $\alpha$ and $X'$ is the closure of the image of $\beta$. 

Let us give a geometric description of $X'$. Let $\Pi\subset X$ be the strict transform of a $2$-plane through the line $\overline{pq}$. The plane $\Pi$ is contracted to a point by the map $\pi_{H_{p,q}}:X\dasharrow \mathbb{P}^{n-2}$ induced by $H_{p,q}$. Indeed, $\pi_{H_{p,q}}$ is induced by the linear projection $\mathbb{P}^n\dasharrow\mathbb{P}^{n-2}$ with center $\overline{pq}$. Observe that a divisor in the linear system of $D'$ has a base component when restricted to $\Pi$; namely the curve $L$. 
\end{tcolorbox}

\begin{tcolorbox}
Therefore, $\phi_{D'|\Pi}$ is the rational map induced by the linear system of conics through $p$ and $q$, hence its image is a smooth quadric surface $Q_{\Pi}\cong\mathbb{P}^1\times\mathbb{P}^1$. The quadric $Q_{\Pi}$ intersects $Z$ at a point. The morphism $\widetilde{\pi}_{H_{p,q}}:X'\dasharrow \mathbb{P}^{n-2}$, induced by the strict transform of $H_{p,q}$ on $X'$, contracts $Q_{\Pi}$ to the point $\pi_{H_{p,q}}(\Pi)$ and maps $Z$ isomorphically onto $\mathbb{P}^{n-2}$. We have the following commutative diagram 
\[
  \begin{tikzpicture}[xscale=1.5,yscale=-1.2]
    \node (A0_0) at (0, 0) {$X$};
    \node (A0_2) at (2, 0) {$X'$};
    \node (A1_1) at (1, 1) {$\mathbb{P}^{n-2}$};
    \path (A0_0) edge [->,swap,dashed]node [auto] {$\scriptstyle{\pi_{H_{p,q}}}$} (A1_1);
    \path (A0_2) edge [->]node [auto] {$\scriptstyle{\widetilde{\pi}_{H_{p,q}}}$} (A1_1);
    \path (A0_0) edge [->,dashed]node [auto] {$\scriptstyle{\varphi_{D'}}$}  (A0_2);
  \end{tikzpicture}
  \]    
and $X'$ has a structure of $(\mathbb{P}^1\times\mathbb{P}^1)$-bundle over $\mathbb{P}^{n-2}$. 

Summing up, the birational model of $X$ corresponding to the chamber $\left\langle H_p,H_q,H_{p,q}\right\rangle$ is a quadric bundle over $\mathbb{P}^{n-2}$ and, as we already noticed, the other chambers $\left\langle H,H_p,E_q\right\rangle$, $\left\langle H,H_q,E_p\right\rangle$ and $\left\langle H,E_p,E_q\right\rangle$ corresponds respectively to $\mathbb{P}^3$ blown-up at $q$, $\mathbb{P}^3$ blown-up at $p$ and $\mathbb{P}^3$. The chamber $\left\langle H,H_p,H_q\right\rangle$ corresponds to $X$ itself.  
\end{tcolorbox}

\section*{Birational twin varieties}\label{bir_twin}
Our aim in this section is to study the weak Lefschetz principle in the context of birational geometry. In other words, if we consider two varieties $X,Y$ equipped with an embedding $i: Y\hookrightarrow X$, then we ask about the objects encoding birational information of $Y$ which are fully determined by $X$. Such objects may include effective and nef cones, or finer information such as the Mori chamber decomposition. The following definition focuses on two possible forms of the weak Lefschetz principle in this context. They were introduced and explored in \cite{LM18}, \cite{LMR18}.

\begin{Definition}\label{BT}
Let $X,Y$ be $\mathbb{Q}$-factorial projective varieties and $i:Y\hookrightarrow X$ be an embedding. These varieties are said to be \textit{Lefschetz divisorially equivalent} if the pull-back $i^*:\Pic(X)\rightarrow\Pic(Y)$ induces an isomorphism such that
$$i^*\Eff(X) = \Eff(Y),\quad  i^*\Mov(X) = \Mov(Y),$$ $$i^*\Nef(X) = \Nef(Y).$$

\noindent
We say that $X$ and $Y$ are \textit{birational twins} if they are Lefschetz divisorially equivalent Mori dream spaces and in addition $i^{*}\MCD(X) = \MCD(Y)$.
\end{Definition}

\begin{Example}
Consider a linear subspace $\mathbb{P}^k\subset\mathbb{P}^n$ passing through the points $p,q$. Let $Y\subset X = Bl_{p,q}\mathbb{P}^n$ be the strict transform of $\mathbb{P}^k$ in $X$. Then $Y$ is isomorphic to the blow-up of $\mathbb{P}^k$ at two points, and via the embedding $i:Y\hookrightarrow X$ all the conditions of Definition \ref{BT} are satisfied when $n,k>1$. In other words, $X$ and $Y$ are birational twins. 
\end{Example}

\noindent
The Grothendieck-Lefschetz theorem 
implies that a smooth variety $X$ of Picard rank one and dimension at least four is Lefschetz divisorially equivalent to any of its effective ample divisors. If the rank of the Picard group is higher than $1$, then B. Hassett, H-W. Lin and C-L. Wang \cite{HLW01} exhibited an example of a variety $X$ with a divisor $D$ such that the inclusion $D \hookrightarrow X$ induces an isomorphism of Picard groups but does not preserve the nef cone; hence, $D$ and $X$ are not Lefschetz divisorially equivalent. In other words, the weak Lefschetz principle fails for the Nef cone in this example.

\noindent
In general, we do not know a classification of the varieties that admit a subvariety which is birational twin to it nor do we know of natural conditions ensuring that varieties $X\subseteq Y$ are birational twins. There are examples where the natural choice fails. However, in the next section we exhibit a series of examples for which the ample cone, the nef cone and even the Mori chamber decomposition satisfy the weak Lefschetz principle.

\section*{Examples: complete quadrics and collineations}
Recently in \cite{LH14, Ma18a, Ma18b, LM18} the birational geometry of classical spaces, called spaces of complete forms, have been studied from the perspective of Definition \ref{BT}. We finish up the article by mentioning some of the ingredients of such a study.

\noindent
Let $V$ be a $K$-vector space of dimension $n+1$ over an algebraically closed field $K$ of characteristic zero. We will denote by $\mathcal{X}(n)$ and $\mathcal{Q}(n)$ the spaces of complete collineations and complete quadrics of $V$, respectively. These spaces are very particular compactifications of the spaces of full rank linear endomorphisms and full-rank symmetric linear endomorphisms of $V$, respectively.

\noindent
In \cite{Va82}, \cite{Va84}, I. Vainsencher showed that these spaces can be understood as sequences of blow-ups of the projective spaces parametrizing $(n+1)\times (n+1)$ matrices modulo scalars along the subvariety parametrizing rank one matrices and the strict transforms of their secant varieties in order of increasing dimension. 

\noindent
Recall that given an irreducible and reduced non-degenerate variety $X\subset\P^N$, and a positive integer $h\leq N$, the \textit{$h$-secant variety} $\sec_h(X)$ of $X$ is the subvariety of $\P^N$ obtained as the closure of the union of all $(h-1)$-planes spanned by $h$ general points of $X$. Spaces of matrices and symmetric matrices admit a natural stratification dictated by the rank. Observe that a general point of the $h$-secant variety of a Segre, or a Veronese, corresponds to a matrix of rank $h$. More precisely, let $\mathbb{P}^N$ be the projective space parametrizing $(n+1)\times (n+1)$ matrices modulo scalars, $\mathbb{P}^{N_{+}}$ the subspace of symmetric matrices, $\mathcal{S}\subset\mathbb{P}^N$ the Segre variety, and $\mathcal{V}\subset\mathbb{P}^{N_{+}}$ the Veronese variety. Since $\sec_h(\mathcal{V}) = \sec_h(\mathcal{S})\cap\mathbb{P}^{N_{+}}$, the natural inclusion $\mathbb{P}^{N_{+}}\hookrightarrow\mathbb{P}^{N}$ lifts to an embedding $i:\mathcal{Q}(n)\hookrightarrow\mathcal{X}(n)$.  

\noindent
We finish up by citing our contribution to the vast realm of Lefschetz-type theorems: the spaces of complete quadrics and complete collineations are Lefschetz divisorially equivalent, for all $n>1$, via the embedding $i:\mathcal{Q}(n)\hookrightarrow \mathcal{X}(n)$ \cite[Theorem A]{LM18}. Furthermore, in the case of complete quadric surfaces, we have that $\mathcal{Q}(3)$ and $\mathcal{X}(3)$ are birational twins \cite[Theorem B]{LM18}.

\begin{Example}
The space $\mathcal{X}(3)$ is the blow-up of the projective space $\mathbb{P}^{15}$ along the Segre variety $\mathcal{S}\cong\mathbb{P}^3\times\mathbb{P}^3\subset \mathbb{P}^{15}$, and also along the strict transform of $\sec_2(\mathcal{S})$. We will denote by $H$ the strict transform of a general hyperplane of $\mathbb{P}^{15}$ and by $E_1,E_2$ the exceptional divisors over $\mathcal{S}$ and $\sec_2(\mathcal{S})$ respectively. Similarly, $\mathcal{Q}(3)$ is the blow-up of the projective space $\mathbb{P}^{9}$ along the Veronese variety $\mathcal{V}\subset \mathbb{P}^{9}$, and also along the strict transform of $\sec_2(\mathcal{V})$. We will denote by $H^{+},E_{1}^{+},E_{2}^{+}$ the divisors on $\mathcal{Q}(3)$ corresponding to the strict transform of a general hyperplane of $\mathbb{P}^9$ and the exceptional divisors over $\mathcal{V}$ and $\sec_2(\mathcal{V})$ respectively. The Mori chamber decomposition $\MCD(\mathcal{X}(3))$ is displayed in the following two dimensional cross-section of $\Eff(\mathcal{X}(3))$:
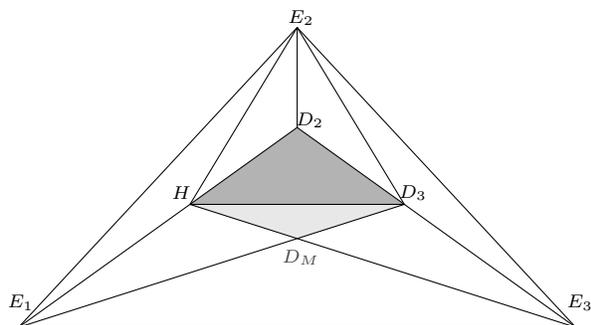
$$
\begin{tikzpicture}[xscale=0.3,yscale=0.7,scale=0.95][line cap=round,line join=round,>=triangle 45,x=1cm,y=1cm]\clip(-13.6,-0.23) rectangle (13.9,6.5);\fill[line width=0pt,fill=black,fill opacity=0.3] (-5.000432432432432,2.4614054054054053) -- (5.000432432432432,2.4614054054054053) -- (0,4) -- cycle;\fill[line width=0pt,color=wwwwww,fill=wwwwww,fill opacity=0.15] (-5.000432432432432,2.4614054054054053) -- (5.000432432432432,2.4614054054054053) -- (0,1.7776389756402244) -- cycle;\draw [line width=0.4pt] (-13,0)-- (13,0);\draw [line width=0.4pt] (13,0)-- (0,6);\draw [line width=0.4pt] (0,6)-- (-13,0);\draw [line width=0.4pt] (0,4)-- (-13,0);\draw [line width=0.4pt] (0,4)-- (13,0);\draw [line width=0.4pt] (-5.000432432432432,2.4614054054054053)-- (13,0);\draw [line width=0.4pt] (5.000432432432432,2.4614054054054053)-- (-13,0);\draw [line width=0.4pt] (-5.000432432432432,2.4614054054054053)-- (5.000432432432432,2.4614054054054053);\draw [line width=0.4pt] (-5.000432432432432,2.4614054054054053)-- (0,6);\draw [line width=0.4pt] (0,6)-- (5.000432432432432,2.4614054054054053);\draw [line width=0.4pt] (0,4)-- (0,6);\begin{scriptsize}\draw [fill=black] (-13,0) circle (0pt);\draw[color=black] (-12.9,0.5) node {$E_1$};\draw [fill=black] (13,0) circle (0pt);\draw[color=black] (13.2,0.5) node {$E_3$};\draw [fill=black] (0,6) circle (0pt);\draw[color=black] (0.18536585365853658,6.2) node {$E_2$};\draw [fill=black] (0,4) circle (0pt);\draw[color=black] (0.6,4.143836565096953) node {$D_2$};\draw [fill=black] (-5.000432432432432,2.4614054054054053) circle (0pt);\draw[color=black] (-5.4,2.7) node {$H$};\draw [fill=black] (5.000432432432432,2.4614054054054053) circle (0pt);\draw[color=black] (5.4,2.7) node {$D_3$};\draw [fill=uuuuuu] (0,1.7776389756402244) circle (0pt);\draw[color=uuuuuu] (0.18536585365853658,1.4) node {$D_M$};\end{scriptsize}\end{tikzpicture}
$$
where $D_M \equiv 6H-3E_1-2E_2$, $D_2\equiv 2H-E_1$, $D_3\equiv 3H-2E_1-E_2$. Here, $E_3\equiv 4H-3E_1-2E_2$ is the class of the strict transform of $\sec_3(\mathcal{S})$. The movable cone $\Mov(\mathcal{X}(3))$ is generated by $\langle H,D_2,D_3,D_M\rangle$. 

\noindent
The spaces $\mathcal{X}(3)$ and $\mathcal{Q}(3)$ are birational twins. Hence, the Mori chamber decomposition $\MCD(\mathcal{Q}(3))$ is obtained from $\MCD(\mathcal{X}(3))$ above by simply replacing $H,E_1,E_2$ with $H^{+},E_1^{+},E_2^{+}$.
\end{Example}

\section*{Acknowledgements}
We thank the referees for their careful reading of our manuscript. Their many constructive suggestions helped us shape this article. 
We also thank Omar Antol\'in Camarena for helpful discussions. During the preparation of this article the first author was partially supported by the CONACYT grant CB-2015/253061.

\bibliography{submitted200915.bib}

\end{document}